\documentclass{article}
\usepackage{amsmath,amssymb,amsfonts}

\def\AAA{{\cal A}}

\newtheorem{theorem}{Theorem}
\newtheorem{conjecture}{Conjecture}

\begin{document}

\title{On the volume of a certain polytope}
\author{Clara S. Chan \and David P. Robbins \and David S. Yuen}
\maketitle

\begin{abstract}
Let $n \ge 2$ be an integer and consider the 
set $T_n$ of $n \times n$ permutation matrices $\pi$ for
which $\pi_{ij}=0$ for $j\ge i+2$.

In this paper we study the convex hull of $T_n$, which we denote by $P_n$.
$P_n$ is a polytope of dimension $\binom{n}{2}$. Our main purpose is
to provide evidence for the following conjecture concerning its volume.
Let $v_n$ denote the minimum volume of a simplex with vertices in the
affine lattice spanned by $T_n$.
Then the volume of $P_n$ is $v_n$ times the product
$$\prod_{i=0}^{n-2} \frac{1}{i+1}\binom{2i}{i} $$
of the first $n-1$ Catalan numbers.

We also give a related result on the Ehrhart polynomial of $P_n$.

\end{abstract}

\section{Introduction}
\label{sec:intro}

Let $n \ge 2$ be an integer and consider the 
set $T_n$ of $n \times n$ permutation matrices $\pi$ for
which $\pi_{ij}=0$ for $j\ge i+2$ and $P_n$ the convex
hull of $T_n$.

Let $V_n$ be the relative volume of $P_n$.  That is, the
volume of $P_n$ expressed in units of the minimum volume $v_n$ of a
simplex with vertices in the affine lattice spanned by $T_n$.
Our main purpose in this paper is to provide evidence for the following
conjecture.

\begin{conjecture}\label{conj:1}
The relative volume $V_n$ of $P_n$ is given by 
$$V_n=  \prod_{i=0}^{n-2} \frac{1}{i+1}\binom{2i}{i}, $$
the product of the first $n-1$ Catalan numbers.
\end{conjecture}

This conjecture arose from the study \cite{CR} of the polytope $B_n$
of all doubly stochastic matrices, which is the convex hull of the set
of all $n \times n$ permutation matrices.
It is easily shown that $P_n$ is a face
of $B_n$ of dimension $\binom{n}{2}$ with $2^{n-1}$ vertices.  
In \cite{CR} we discuss two methods for finding the volume of $B_n$ and
its faces.  We assume some familiarity with these methods, which apply
to the calculation of the volume of $P_n$.
The reader may also wish to consult the other references listed here, which
were central to the work in \cite{CR}.

The first method discussed in \cite{CR} consists of decomposing the polytope
into simplices, each of volume $v_n$, and counting the simplices. By adapting
the method slightly we were able easily to find the relative volumes of $P_n$
and its faces provided that $n \le 10$.
This provided the first evidence for Conjecture \ref{conj:1}.

The second method discussed in \cite{CR} computes the Ehrhart polynomial of
the polytope.  In general
the Ehrhart polynomial of a $d$-dimensional polytope $P$ with integer
vertices is a degree $d$ polynomial (in $t$) denoted $e(P,t)$, with
the property that the number of integer points in the polytope $t\cdot
P$ is $e(P,t)$ when $t\ge 0$.  A basic property of the Ehrhart
polynomial is that the relative volume of the polytope is given by
$d!$ times its leading coefficient.  A
common method for computing the Ehrhart polynomial is to count the
numbers of lattice points in $t\cdot P$ for small $t$ and then to find
the polynomial by interpolation.  For a typical face of $B_n$ 
the Ehrhart polynomial method seems to be more expensive than the
simplicial decomposition method.  However for $B_n$ itself
the Ehrhart polynomial method is less expensive because it is
possible to exploit the symmetries of $B_n$.  These symmetries do not
help with the calculation of the Ehrhart polynomial of $P_n$.  However,
different simplifications in the case of $P_n$ allow us to compute the Ehrhart
polynomial and thus verify Conjecture \ref{conj:1} for $n \le 12$, as described
in Section 2 of this paper.

In Section 3 of this paper we give a proof of a bijection between the
simplices in a decomposition of $P_n$ and a set of easily described
integer arrays, which suggest that a combinatorial proof of Conjecture
\ref{conj:1} may exist.  We also discuss a generalization of the conjecture
which arises from the bijection.

In Section 4 we discuss formulas for the relative volumes of some of the
facets of $P_n$, which we pursued as an alternative path toward proving
Conjecture \ref{conj:1}.  The formulas were discovered by using the
simplicial decomposition method mentioned above.

\section{The Ehrhart polynomial of $P_n$}

One approach to calculating the volume of $P_n$ is to calculate
its Ehrhart polynomial.  Let us denote the Ehrhart polynomial 
of $P_n$ (evaluated at $t$) by $e(P_n,t)$.  Then $e(P_n,t)$ is
the number of ways of filling a left-justified array of $n$ rows
of lengths  $2,3,\ldots,n-1,n,n$ with non-negative integers in
such a way that all row and column sums are $t$.  Thus $e(P_3,1)=4$
since the only four suitable arrays are
\[
\begin{array}{ccc}
   1 & 0  \\
   0 & 1 & 0 \\
   0 & 0 & 1 \\
\end{array}
\quad 
\begin{array}{ccc}
   1 & 0  \\
   0 & 0 & 1 \\
   0 & 1 & 0 \\
\end{array}
\quad
\begin{array}{ccc}
   0 & 1  \\
   1 & 0 & 0 \\
   0 & 0 & 1 \\
\end{array}
\quad
\begin{array}{ccc}
   0 & 1  \\
   0 & 0 & 1 \\
   1 & 0 & 0 \\
\end{array}
\]
$e(P_n,t)$ is known be a polynomial in $t$ whose degree as
a function of $t$ is the dimension of $P_n$ or $\binom{n}{2}$.

We have $e(P_n,0)=1$ for all $n$.  Also it is easily verified
that $e(P_n,1)=2^{n-1}$ for all $n$.  These are special cases
of a more general principal.  

\begin{theorem}\label{thm:1}
For every non-negative integer $t$, the sequence 
\[
e(P_0,t),e(P_1,t),\ldots,e(P_n,t),\ldots
\]
satisfies a linear recursion of degree
$p(t)$ with integer coefficients, where $p(t)$ is the number of
partitions of $t$.
\end{theorem}

We remark that it is conceivable, as far as we know, that the
sequences also satisfy recursions of lower degree.  However, for 
all cases $t=0,\ldots,12$, where we have computed the coefficients of
the linear recursion, the associated characteristic polynomial has
been irreducible over the integers so, in these cases, no lower degree
recursion exists.

\medskip
\noindent
{\em Proof:}  Let us fix a non-negative integer $t$.  Let 
$\pi=(x_1,x_2,\ldots,x_l)$ be a partition of $t$ of length $l \ge 1$.
That is, $0 < x_0\le x_1 \le \cdots \le x_l$ and $x_1+\cdots+x_l=t$.
For integers $n \ge 2$ let $F(\pi,n)$ denote the set
of arrays of $n$ left-justified rows of non-negative integers of 
lengths $l+1,\ldots,l+n-1,l+n-1$ such that the first $l$ column
sums are $x_1,\ldots,x_l$, the remaining column sums are $t$, and
all row sums are $t$.  We set $f(\pi,n)$ to be the cardinality of 
$F(\pi,n)$, and $f(\pi,1)=1$.  Note that when
$\pi$ is the one-part partition $(t)$, we have $f(\pi,n)=e(P_n,t)$.

Suppose that $n \ge 2$.  Set $x_{l+1}=t$.  Let $y_1,\ldots,y_{l+1}$ be any
non-negative integers with $y_i \le x_i$ for $i=1,\ldots,l+1$ such that
$y_1+\cdots+y_{l+1}=t$.  
If $ y_1,\ldots,y_{l+1}  $ is the first row of one of the arrays of $F(\pi,n)$,
then the rest of the array has its first $l+1$ column sums equal to
$z_i=x_i-y_i$.  By deleting the $z_i$'s which equal 0 and sorting the
remaining $z_i$'s, we obtain another partition $\sigma$ of $t$, of length
at most $l+1$, and the number of ways of completing the array is clearly
$ f(\sigma,n-1)$.  (Since $f(\sigma,1)=1$, this also holds for $n=2$.)
Now for every partition $\sigma$, let $M(\pi,\sigma)$ denote
the number of $(l+1)$-tuples $y_1,\ldots,y_{l+1}$ for which our process
of forming the $z$'s by subtracting from the $x$'s, deleting 0's, and sorting 
yields the partition $\sigma$.  Then we have shown that
$$  f(\pi,n)=\sum_{\sigma} M(\pi,\sigma)f(\sigma,n-1).  $$
For fixed $n$ we can regard the array of $f(\pi,n)$, as $\pi$ varies over
partitions of $t$, as a column vector of integers of length $p(t)$.
When $n=1$ we have the vector of all 1's.  The preceding equation shows that
the $n^{th}$ vector is obtained by applying the matrix $M^{n-1}$ to
the all 1's vector.  Thus the sequence of vectors satisfies a linear
recursion with integer coefficients given by the characteristic polynomial
of the matrix $M$.  In particular, the component of the column vector
$f(\pi,n)$ corresponding to the partition $(t)$ is $e(P_n,t)$,
which proves our theorem.  

\medskip
\noindent
{\em Example:}
It is easy to compute the matrix $M$ for small values of $t$.
For example, let $t=2$.  If $\pi=(2)$, then $(x_1,x_2)=(2,2)$,
so $(y_1,y_2)=(0,2)$, $(1,1)$ and $(2,0)$, 
which yield $\sigma=(2)$, $(1,1)$ and $(2)$ respectively.
If $\pi=(1,1)$, then $(x_1,x_2,x_3)=(1,1,2)$, so $(y_1,y_2,y_3)=(0,0,2)$,
$(0,1,1)$, $(1,0,1)$, and $(1,1,0)$, which yield 
$\sigma=(1,1)$, $(1,1)$, $(1,1)$, and $(2)$, respectively.  Thus we have
\[
M=\left(
\begin{array}{rr}
     2  & 1 \\
     1  & 3 \\
\end{array}
\right),
\]
with the rows and columns of $M$ indexed by the partitions
$(2)$ and $(11)$ in that order.
The characteristic polynomial of $M$ is $\lambda^2-5\lambda+5$.  Thus
we have  $e(P_n,2)=5e(P_{n-1},2)-5e(P_{n-2},2)$.  Initial values are
$e(P_1,2)=1$ (by definition) and $e(P_2,2)=3$ (by applying $M$ to the
all ones vector).

\medskip

Theorem \ref{thm:1} and its example contain the essential ideas behind
our method for evaluating $e(P_n,t)$ for small values of $t$.
If we wish to calculate a value of $e(P_n,t)$ for which $n$ is also small,
we can simplify a little more by computing and using only the submatrix of
$M$ corresponding to partitions of length not exceeding $n$.

Let us denote the characteristic polynomial of the matrix associated to the
nonnegative integer $t$ by $f_t(\lambda)$.  Here are the first 6 polynomials.
\begin{eqnarray*}
 f_0&=&\lambda-1  \\
 f_1&=&\lambda-2  \\
 f_2&=&\lambda^2-5\lambda+5 \\
 f_3&=&\lambda^3-10\lambda^2+27\lambda-20 \\
 f_4&=&\lambda^5-20\lambda^4+135\lambda^3-396\lambda^2+518\lambda-245 \\
 f_5&=&\lambda^7-36\lambda^6+480\lambda^5-3140\lambda^4+11059\lambda^3-21180\lambda^2+20560\lambda-7840 \\
\end{eqnarray*}
As far as we have computed, all the roots of these polynomials are positive
real numbers.

We can also use Ehrhart's reciprocity principle to simplify the
computation of the Ehrhart polynomial.  Recall that,
for a $d$-dimensional polytope $P$ with integer vertices, and
$t>0$, Ehrhart's reciprocity principle states that
\[
e^*(P,t)=(-1)^d e(P,-t)
\]
where $e^*(P,t)$ is the the number of lattice points in the 
interior of $t\cdot P$. 

An interior lattice point of $t\cdot P_n$  is an array of positive integers
consisting of left-justified rows of length $2,3,\ldots,n,n$ with all row and
column sums equal to $t$.  For such an array, if $k<n$, the first $k$ rows
have lengths $2,3,...,k,k+1$ and the sum of all their entries taken together
is $tk$.  On the other hand, the sum of all entries in the first $(k+1)$
columns is $(k+1)t$, and this includes all entries in the first $k$ rows.
Thus the sum of all entries in the first $k+1$ columns of the last $n-k$ rows
must be $t$.  Since all entries are positive, it follows that 
$ t \ge (k+1)(n-k) $ and this inequality must hold for $k=0,\ldots,n-1$.
Thus if an interior point exists for a given $t$, we must have $t$ at least equal to the maximum over $k$ of $(k+1)(n-k)$.
Thus for odd $n=2m+1$, we have $e(P_n,-t)=0$ for $t=1,\ldots,(m+1)^2-1$,
while for even $n=2m$, we have $e(P_n,-t)=0$ for $t=1,\ldots,m(m+1)-1$.

We have calculated the Ehrhart polynomials of $P_n$ for $n=2,\ldots,12$.
Here are the first few.
\begin{eqnarray*}
e(P_2,t)&=&t+1 \\
e(P_3,t)&=&\frac{1}{6}\prod_{i=1}^3 (t+i) \\
e(P_4,t)&=&\frac{t+3}{360}\prod_{i=1}^5(t+i) \\
e(P_5,t)&=&\frac{(t+3)^2}{362880}\prod_{i=1}^8(t+i) \\
e(P_6,t)&=&\frac{(t+3)^2(t^2+12t+26)}{9340531200} \prod_{i=1}^{11}(t+i) \\
e(P_7,t)&=&\frac{(t+3)^2(14t^4+353t^3+2985t^2+9568t+10336)}{121645100408832000}\prod_{i=1}^{15}(t+i)\\
\end{eqnarray*}
The factor $(t+3)^2$ which appears in $e(P_5,t)$ persists through
$n=12$, but we have no proof that it persists forever.

To check Conjecture \ref{conj:1}, one multiplies the leading coefficient of
$e(P_n,t)$ by $\binom{n}{2}!$, to get the predicted relative volume.
This works through $n=12$.

\section{Explicit Decomposition into simplices}

In this section we show that the polytope $P_n$ can be decomposed into 
minimal volume simplices which are in bijection with an easily described set of integer
arrays.  Thus the relative volume of $P_n$ is simply the number of such
integer arrays.  This suggests an avenue for proving Conjecture \ref{conj:1},
although we have not been successful thus far.  Postnikov and Stanley \cite{PS}
have found a bijection very much like ours, and also observed that therefore
Conjecture \ref{conj:1} is 
equivalent to
$$ K\left(a_1 + 3a_2 + 6a_3 + \cdots + \binom{n}{2} a_{n-1}\right) = \prod_{i=0}^{n-1} \frac{1}{i+1}\binom{2i}{i} $$
where $a_1,\ldots,a_{n-1}$ 
is a choice of simple roots and $K$ is the Kostant partition function 
for the root system $A_{n-1}$.

To describe our decomposition we first need some notation.

The polytope $P_n$ consists of doubly stochastic matrices
$Y=(y_{ij})$ where $y_{ij}=0$ for $j> i+1$.

However the entries $y_{ij}$, $j\le i \le n-1$, determine
the remaining $2n-1$ entries 
$y_{12},y_{23},\ldots,y_{n-1,n}$
and $y_{n1},\ldots,y_{nn}$.  Thus
we may view a point $Y$ in $P_n$ as a triangular array
$$
\begin{array}{cccccc}
y_{11} &  &  &  &  & \\
y_{21} & y_{22} &  &  &  & \\
\vdots &        &        &        &   & \\
y_{n-1,1} & y_{n-1,2} & y_{n-1,3} & \cdots & y_{n-1,n-1} & 
\end{array}.
$$
where the nonnegative $y_{ij}$'s satisfy the conditions

\begin{equation} \label{eq:1}
\sum_{i=k}^{n-1}y_{ik}\le \sum_{j=1}^{k-1}y_{k-1,j}\le 1
\end{equation}
for $k=2,\ldots,n-1$ and 
\begin{equation} \label{eq:2}
\sum_{i=1}^{n-1}y_{i1}\le 1
\end{equation} so that the first column has sum $\le 1$.

Let $\AAA_n$ be the set of triangular arrays of non-negative
integers 
$$
\begin{array}{cccccc}
a_{22} &  &  & \\
a_{32} & a_{33} &  & \\
\vdots & \vdots   & \ddots       &        & \\
a_{n-1,2} & a_{n-1,3} & \cdots & a_{n-1,n-1}
\end{array}
$$
where the $a$'s are subject to the constraint that
\begin{eqnarray} 
 a_{22}+\cdots+a_{n-1,2} &\le& 0 \nonumber \\
 a_{33}+\cdots+a_{n-1,3} &\le& 1+a_{22} \nonumber  \\   
 a_{44}+\cdots+a_{n-1,4} &\le& 2+a_{32}+a_{33}    \label{eq:3x} \\
 a_{55}+\cdots+a_{n-1,5} &\le& 3+a_{42}+a_{43}+a_{44} \nonumber \\
  \vdots                                 \nonumber     \\
 a_{n-1,n-1}             &\le& n-3+ a_{n-2,2}+a_{n-2,3}+\cdots+a_{n-2,n-2} \nonumber 
\end{eqnarray}
Note that this condition implies that the leftmost column in the
array is all zeros.

For example $\AAA_5$ consists of the 10 triangular arrays
$$
\begin{array}{cccc}

\begin{array}{ccc}
0 & & \\
0 & 0 & \\
0 & 0 & 0 \\
\end{array}
&
\begin{array}{ccc}
0 & & \\
0 & 0 & \\
0 & 0 & 1 \\
\end{array}
&
\begin{array}{ccc}
0 & & \\
0 & 0 & \\
0 & 0 & 2 \\
\end{array}
&
\\
\phantom 1 & & & \\
\begin{array}{ccc}
0 & & \\
0 & 0 & \\
0 & 1 & 0 \\
\end{array}
&
\begin{array}{ccc}
0 & & \\
0 & 0 & \\
0 & 1 & 1 \\
\end{array}
&
\begin{array}{ccc}
0 & & \\
0 & 0 & \\
0 & 1 & 2 \\
\end{array}
&
\\
\phantom 1 & & & \\
\begin{array}{ccc}
0 & & \\
0 & 1 & \\
0 & 0 & 0 \\
\end{array}
&
\begin{array}{ccc}
0 & & \\
0 & 1 & \\
0 & 0 & 1 \\
\end{array}
&
\begin{array}{ccc}
0 & & \\
0 & 1 & \\
0 & 0 & 2 \\
\end{array}
&
\begin{array}{ccc}
0 & & \\
0 & 1 & \\
0 & 0 & 3 \\
\end{array}
\\

\end{array}
$$

We will give a decomposition of $P_n$ into simplices 
all of the same volume
in such a way that the simplices in the decomposition will
be in one-to-one correspondence with the set $\AAA_n$. 

We start by defining a mapping which assigns to
each element $\alpha$ of $\AAA_n$ a simplex contained in $P_n$.

Let $V$ be the space of triangular arrays
$$
\begin{array}{cccccc}
x_{11} &  &  &  &  & \\
x_{21} & x_{22} &  &  &  & \\
\vdots & \vdots &   \ddots     &  &   & \\
x_{n-1,1} & x_{n-1,2} & \cdots & x_{n-1,n-1} & 
\end{array}.
$$
Let $S$ be the unit simplex in $V$ consisting of all nonnegative
arrays of the preceding form in which the sum of all the entries is 
$\le 1$.  The simplex associated to $\alpha$ will be the image of
$S$ under a certain linear transformation of determinant 1 which is 
associated to
$\alpha$.  We will denote this linear transformation by $L(\alpha)$.  

We can construct many such mappings from $\alpha$'s to
simplices, each of which yields a suitable
decomposition of $P_n$.  A convenient way
to specify a single one of these is to assign a linear
ordering to the variables $x_{ij}$.  It does not
matter what linear ordering we use.

To form $L(\alpha)$ we start with the identity matrix, represented
by the preceding triangle.  Then $L(\alpha)$ is formed by a series
of steps, one step for each column of $\alpha$.  At the beginning of
each step we have a linear transformation consisting of a triangle
of linear functions of the $x$'s.  The step itself consists in 
performing certain operations on the triangle, leading to
another triangle of linear functions.

After each step all the linear functions in each triangle are
of a particularly simple form in that

(C0) 
each linear function is a linear combination
of the $x$'s all of whose coefficients are 0 or 1; i.e.
a sum of distinct $x$'s.  Moreover no two entries in any row involve
the same variable and no two entries in any column involve the same
variable.

Also, after having used columns $2,\ldots,k$ of $\alpha$, the triangle 
of linear functions has certain additional properties depending on 
$k$. 

(C1)
In the rectangular sub-array consisting of columns $1,\ldots,k$ of
rows $k$ through $n-1$, no two of the linear functions share any variables.  

(C2)
Within
this rectangle, if $j \ge 2$,
entry $ij$ is a sum of precisely $a_{ij}+1$ variables, while entry $i1$ is
just $x_{i1}$, a sum of 1 variable. 

(C3)
The only variables that occur in columns $1,\ldots,k$ of the 
triangle are those that were originally in these columns.

(C4)
In columns $k+1,\ldots,n-1$, the linear function is just
the original variable $x_{ij}$.

(C5) For $2 \le j \le k$ the variables appearing in column $j$ of the array
are a proper subset of those appearing in row $j-1$.

(C6)
In the first column of the triangular array every
variable originally in the first $k$ columns
appears precisely once.

In view of (2) above, after having used columns $2,\ldots,k$ of
$\alpha$, 
there are precisely $k+a_{k2}+\cdots+a_{k,k}$ variables
in all the sums in the $k$th row.  Let us denote this number of variables by $N$.
We now list these variables in the assigned order, denoting them as
$z_1,\ldots,z_N$.

From our conditions above defining $\AAA_n$, we have 
\begin{equation}\label{eq:3}
N > a_{k+1,k+1}+\cdots+a_{n-1,k+1}
\end{equation}

Before actually modifying the triangle of linear functions we
first use column $k+1$ of $\alpha$ to 
parse the $z$'s into ``chunks'', putting the first
$ a_{k+1,k+1}$ $z$'s into the first chunk, the
next $a_{k+2,k+1}$ $z$'s into the next chunk, and so
forth, one chunk for each entry
$a_{i,k+1}$ in column $k+1$ of $\alpha$.  The inequality (\ref{eq:3x}) guarantees that
there is at least one more variable available than is needed to form all the
chunks.
Notice that some of
the chunks can be empty and that at least one of the
$z$'s does not appear in any chunk.  

After we form each 
chunk, we associate to it
the first of the $z$'s (in our
ordering) that has not yet appeared in any chunk (including the chunk
just formed) and call this
the ``cap'' of the chunk just formed.  Note that since some of the
chunks are empty, it is possible for chunks associated to 
several consecutive entries $a_{i,k+1}$ to have the same
cap.

Now we modify the triangle of linear functions in two
substeps.

First, for each $i=k+1,\ldots,n-1$, we let $z$ be the cap of the chunk associated
to $a_{i,k+1}$, and
then replace every occurrence of $z$ in columns
$1,\ldots,k$ of the triangle of linear functions 
with $z+x_{i,k+1}$.  The order in which we perform the substitutions 
in this substep
is immaterial since the variables $x_{i,k+1}$ did not
previously appear in columns $1,\ldots,k$ (and hence also are never
caps).

Second, for each variable
$x_{i,k+1}$ in column $k+1$ of the triangle of linear functions, we
replace that variable by a sum consisting of the
variable itself plus the sum of all the variables 
in the chunk associated to $a_{i,k+1}$.

Each of the conditions (C0)-(C6) holds for the initial
triangle and it is easy to see, inductively, that
the modification rules above preserve the conditions.
Thus they hold at every stage.  

A somewhat deeper property of our inductive procedure is
that, at each stage, the triangle of linear functions represents
a linear transformation of determinant 1.

The first substep is a linear substitution of determinant 1.
However, the second substep is not strictly a linear substitution
since the first substep results in occurrences $x_{i,k+1}$ in
the columns $1,\ldots,k$, while, in the second substep, we do not
perform the substitution in these occurrences.  

However, we can obtain the second substep by a sequence of
pairs of linear substitutions 
of determinant 1.  Indeed for each $x_{i,k+1}$ in
column $k+1$ we substitute for its cap variable the cap variable minus
the sum of the associated chunk variables, and then substitute for
the $x_{i,k+1}$, $x_{i,k+1}$ plus the sum of the chunk variables.  The
effect of the these two substitutions is to leave columns $1,\ldots,k$ unchanged
but to perform the desired substitution in column $k+1$.

Here is an illustration of the procedure described above.
Suppose that $n=5$ and that the array
$\alpha$ is 
$$
\begin{array}{cccc}
0 &   &   &    \\
0 & 1 &   &    \\
0 & 0 & 2 &    \\
0 & 0 & 1 & 2  \\
\end{array}.
$$
Let us take the array of variables (called $x_{ij}$ above)
to be 
$$
\begin{array}{lllll}
A &   &     &   &   \\
B & F &     &   &   \\
C & G & J   &   &   \\
D & H & K   & M &   \\
E & I & L   & N & O \\
\end{array}
$$
and let us order the variables alphabetically.

When the first column of $\alpha$ is used, all 4 chunks have length
0 and cap $A$.  So the effect is that all four variables are added to $A$,
yielding
$$
\begin{array}{lllll}
AFGHI &   &   &   &  \\
B     & F &   &   &  \\
C     & G & J &   &  \\
D     & H & K & M &  \\
E     & I & L & N & O\\
\end{array}.
$$
where, for the rest of this example, we designate addition by
juxtaposition so that we abbreviate $A+F+G+H+I$ with $AFGHI$.

When the second column of $\alpha$ is used, the variables in the
second row of the triangle are $B,F$ and there are three
chunks, the first is $B$ and the last two are empty.  All
three have cap $F$.  Thus we obtain 
$$
\begin{array}{lllll}
AFGHIJKL &      &    &   &   \\
B        & FJKL &    &   &   \\
C        & G    & BJ &   &   \\
D        & H    & K  & M &   \\
E        & I    & L  & N & O \\
\end{array}.
$$

When the third column of $\alpha$ is used there, are four variables
in the third row of our triangle, namely $B,C,G,J$ and there are two chunks,
$B,C$ and $G$, with caps $G$ and $J$ respectively.  The chunks
are adjoined to $M$ and $N$ and, in the first three columns of the
triangle, $G$ is replaced by $GM$ and $J$ by $JN$.  Thus we obtain
$$
\begin{array}{lllll}
AFGHIJKLMN &       &     &     &   \\
B          & FJKLN &     &     &   \\
C          & GM    & BJN &     &   \\
D          & H     & K   & BCM &   \\
E          & I     & L   & GN  & O \\
\end{array}.
$$

Finally, when the last column of $\alpha$ is used, there are 6 variables
in the fourth row of the triangle, $B,C,D,H,K,M$.  We form one
chunk of size 2, namely $B,C$, with cap $D$, obtaining
$$
\begin{array}{lllll}
AFGHIJKLMN &       &     &     &    \\
B          & FJKLN &     &     &    \\
C          & GM    & BJN &     &    \\
DO         & H     & K   & BCM &    \\
E          & I     & L   & GN  & BCO\\
\end{array}.
$$

The unit simplex in $V$ is the set of 15-tuples $A,\ldots,O$ of non-negative
reals whose sum is $\le 1$. 
Still, taking note of our juxtaposition notation for addition, we see that
the triangle above defines a linear mapping from the unit simplex to $P_5$.  

It is easy to see, inductively, that this will be the case for any $\alpha$ in 
$\AAA_n$.  First note
the inequality (\ref{eq:2}) will always hold because of (C6).  One
also easily verifies that the conditions (\ref{eq:1}) always hold.  The
second of the inequalities is a consequence of the fact that the
variables occurring in any row are always distinct.  The first inequality 
is a consequence of (C5).

Thus we have associated to every $\alpha$ in $\AAA_n$ a simplex whose
volume is $1/\binom{n}{2}!$.  

One needs also to show that the simplices $L(\alpha)$ cover
$P_n$ and have disjoint interiors.
There is an argument, rather similar to the preceding, in which we
start with a point of $P_n$ and build up $\alpha$ and $L(\alpha)$
with a construction like the preceding.  But we omit the details.

Thus our conjecture would be proved if we could show that the 
cardinality of $\AAA_n$ was given by
$$\prod_{i=0}^{n-2} \frac{1}{i+1}\binom{2i}{i}. $$
We have not been able to show this. 

However, this combinatorial
interpretation leads to a stronger conjecture.  We can classify 
the elements of $\AAA_n$ according to the number of times that
we have equality in (\ref{eq:3x}).  This can hold from 1 to $n-2$
times.

\begin{conjecture}\label{conj:3a}

Suppose that $n \ge 2$ and that $D_{nk}$ is the number of elements of $\AAA_n$ 
for which equality holds for $k$ of the inequalities (\ref{eq:3x}).
Then $D_{nk}$ is divisible by 
$$\prod_{i=0}^{n-3} \frac{1}{i+1}\binom{2i}{i} $$
and the quotient is the {\em Narayana number} $N(n-2,k)$, where
$$  N(n,k)=\frac{1}{n}\binom{n}{k}\binom{n}{k-1}. $$

\end{conjecture}

For example, the following two elements of $\AAA_5$ satisfy just
1 equality in (\ref{eq:3x})
$$
\begin{array}{ccc}
0 & & \\
0 & 0 & \\
0 & 0 & 0 \\
\end{array}
\quad
\begin{array}{ccc}
0 & & \\
0 & 0 & \\
0 & 0 & 1 \\
\end{array}
$$
while the following two elements of $\AAA_5$ satisfy three
equalities
$$
\begin{array}{ccc}
0 & & \\
0 & 0 & \\
0 & 1 & 2 \\
\end{array}
\quad
\begin{array}{ccc}
0 & & \\
0 & 1 & \\
0 & 0 & 3 \\
\end{array}
$$
The remaining 6 satisfy two equalities.

To conclude this section we define a generalization of the set $\AAA_n$
and a corresponding generalization of Conjecture \ref{conj:1}.

Let $\AAA_n^j$ denote the set of elements of $\AAA_n$ in which the first
$j$ columns consist entirely of zeros.  Thus $\AAA_n^1 = \AAA_n$ and
$\AAA_n^j \subset \AAA_n^{j-1}$ for all $j\ge 2$.
Here is a small table of values of $\AAA_n^j$.

$$
\begin{array}{rrrrrr}
n & j=1 & j=2 & j=3 & j=4 & j=5 \\
3 & 1 & & & & \\
4 & 2 & 1 &  & & \\
5 & 10 & 3 & 1 & &\\
6 & 140 &  28 &  4 &   1 & \\
7 & 5880 & 840 & 60 &  5 & 1 \\
\end{array}
$$

We have the following

\begin{conjecture}\label{conj:3b}
The number of elements in $\AAA_n^j$ is the product

$$\prod_{i=j}^{n-3} \frac{1}{2i+1}\binom{n+i-1}{2i}. $$

\end{conjecture}

With the help of the computer program Mathematica we can verify this 
easily for $n-j \le 6$.

\section{Facets of $P_n$ and their volumes}

Another approach toward proving Conjecture \ref{conj:1} is to try to
understand the relative volumes of the facets of $P_n$.  In this section we
study these facet volumes and make a conjecture concerning these volumes
based on evidence obtained by the simplicial decomposition method described
in \cite{CR}.

Suppose that $n \ge 2$ is an integer and that
$1 \le r,s \le n$ and $s \le r+1$.  
Consider the convex hull $P_n(r,s)$ of those permutations in 
$T_n$ whose $(r,s)$ entry is zero.  Then $P_n(r,s)$ is
always a face of $P_n$.  If $n=2$ these are all facets of
$P_n(r,s)$,  but for $n\ge3$, $P_n(r,s)$ is a facet of $P_n$ precisely
when $r \ne 1$ and $s \ne n$ and $s \ne r+1$.

Since the set $T_n$ is invariant under the operation of 
exchanging the first two columns and the operation of 
exchanging the last two rows, the same symmetries apply to 
the volumes of the facets. 
Thus the volume of $P_n(r,1)$ is equal to that of $P_n(r,2)$ for all $r$.
Also the volume of $P_n(n,s)$ equals that of $P_n(n-1,s)$.  Thus, we can
display the volumes of all the facets as a triangular array consisting of the
volumes of $P_n(r,s)$ for 
$2 \le r \le n-1$ and $2 \le s \le n-1$ and $s \le r$.

Here are the volumes of $P_n(r,s)$ for $ n=3,\ldots,7 $.  
\[
\begin{array}{rrrrr}
\\
       1  \\
\\
       1  \\
       2 &    1 \\
\\
       3 &  \\
       7 &    4  \\
      10 &    7 &     3 \\
\\
      28 \\
      70 &    42 \\
     112 &    84 &   42  \\
     140 &  112  &   70 &    28  \\
\\
     840  \\
    2180 &  1340  \\
    3700 &  2860 &  1520  \\
    5040 &  4200 &  2860 &  1340  \\
    5880 &  5040 &  3700 &  2180  & 840 \\
\end{array}
\]
These arrays have some properties that are easily verified.  For
example, there is symmetry about the anti-diagonal.  
There is a slightly deeper fact.  In any $2 \times 2$ submatrix
of the preceding array the sum of the entries on one diagonal
of the submatrix is equal to the sum of the entries on the
other.

There is a slightly stronger version that can
be stated a little
more elegantly if we add an extra diagonal of  zeroes above the
main diagonal and then complete the triangle to a skew-symmetric
matrix.  For example, the square matrix associated to the last
triangle (corresponding to $n=7$) is 
\[
\begin{array}{rrrrrr}
       0 & -840  & -2180 & -3700 & -5040 & -5880 \\
     840 &    0  & -1340 & -2860 & -4200 & -5040 \\
    2180 &  1340 &     0 & -1520 & -2860 & -3700 \\
    3700 &  2860 &  1520 &     0 & -1340 & -2180 \\
    5040 &  4200 &  2860 &  1340 &     0 &  -840 \\
    5880 &  5040 &  3700 &  2180 &   840 &     0 \\
\end{array}
\]
Each square matrix formed this way has the property
that, for any of its $2 \times 2$ submatrices, the sum on the
entries on one diagonal is the same as the sum of the
entries on the other.

This property is easily proved.  It results from the fact that
the relative volume of $P_n$ can be expressed as the sum of the
relative volumes of the facets opposite any vertex.  The
rectangular relations arise from pairs of vertices that 
(when regarded as permutations) differ by a transposition.

One consequence of the rectangular relations is that all 
the entries in each triangle depend linearly on the main
diagonal so we can describe the whole triangle, much more
succinctly,  in terms of the diagonal.
Here are
the diagonals (listed as rows) for $n=3,\ldots,10$.
(The last four rows need to be
completed to be palindromic of length $n-2$.)
\[
\begin{array}{rrrrr}
 1 \\
 1 & 1\\
 3 & 4 &3  \\
 28 & 42 &42 &28\\
 840 & 1340 & 1520 &1340 \\
 83160 & 137610 & 167310 &167310 \\
 27747720 & 47016970 & 59676120 & 64091020\\
 31743391680 & 54669174560 & 71411118240 & 80251753120 \\
\end{array}
\]
The entries in the first two columns of this array seem to be predictable.  Suppose
that $a_n$ denotes the entry in the first column and $b_n$ the 
entry in the second column.  Then
$a_n$ is defined for $n\ge 3$, and $b_n$ is defined for $n \ge 4$ 
so that 
$$a_3,a_4,a_5,a_6\ldots=1,1,3,28 \ldots$$
 and 
$$b_4,b_5,b_6,b_7\ldots=1,4,42,\ldots. $$

\begin{conjecture}\label{conj:2a}
For $n \ge 3$,
$$    a_n=3V_n/\binom{n}{2}.  $$
For $n \ge 4$,
$$   (n-1)\left( \frac{b_{n+1}}{a_{n+1}}-\frac{b_n}{a_n}\right) = 
(n+2)\left(\frac{b_{n+2}}{a_{n+2}}- \frac{b_{n+1}}{a_{n+1}}  \right).
$$
\end{conjecture}

The evidence for the second formula is perhaps not all that
compelling since the result is known to hold only for $n=4,\ldots,8$.
However it is not hard to check that the two formulas above, taken together
with Conjecture \ref{conj:1}, predict integral values for
$b_n$, for all $n$.  So this gives some additional evidence.

\end{document}